\documentclass[12pt]{article}
\usepackage{graphicx}
\usepackage{amssymb}
\usepackage{amsmath}
\numberwithin{equation}{section}

\newtheorem{Pa}{Paper}[section]
\newtheorem{Tm}[Pa]{{\bf Theorem}}
\newtheorem{La}[Pa]{{\bf Lemma}}
\newtheorem{Cy}[Pa]{{\bf Corollary}}
\newtheorem{Rk}[Pa]{{\bf Remark}}

\newtheorem{Dn}[Pa]{{\bf Definition}}
\newtheorem{Pn}[Pa]{{\bf Proposition}}

\title{Eigenfunction expansions and scattering theory associated with Dirac equation}
\author{Lev Sakhnovich}

\date{}

\begin{document}
\maketitle

\emph{99 Cove ave. Milford, CT,06461, USA}

\vspace{0.3em}
 
 E-mail: lsakhnovich@gmail.com
 
 \vspace{0.3em}
 
\begin{abstract}The classical Lippmann-Schwinger equation (LS equation) plays an important role in the
scattering theory (Schr\"odinger   equation, non-relativistic case). In  our previous paper \cite{Sakh4} we considered the relativistic analogue of the Lippmann-Schwinger equation (RLS equation). We presented the RLS equation
in the integral form. In the present paper, we use the   RLS equation in order to study the scattering problems (stationary and dynamical cases)
for relativistic Dirac equation.
Our approach allows us  to develop the RLS equation  theory that is comparable in its completeness with  the LS equation theory. In particular, we consider the eigenfunction expansion associated with the relativistic Dirac equation.
 We note that the works on the theory of the LS equation (see \cite{IK, RS, Sim}) serve as a model for us.\end{abstract}
 
  

\section{Introduction}The classical integral
Lippmann-Schwinger equation (LS-equation) plays an important role in the  scattering theory (non-relativistic case,  Schr\"odinger   equation).
 The relativistic analogue of the Lippmann-Schwinger equation  was formulated in the terms of the limit
values of the corresponding resolvent. In the our previous paper \cite{Sakh4} we found the limit values of the resolvent
in the explicit form. Using this result we have represented relativistic Lippmann-Schwinger equation (RLS equation) as an integral equation.
In the present paper we use the corresponding integral equation  and investigate the scattering problems, stationary and dynamical cases.
 Our approach allows  to develop a RLS equation  theory that is comparable in its completeness to the theory of LS equation. In particular we consider the eigenfunction expansion associated with relativistic Dirac equation. We note that the works on the theory of the LS equations Ikebe \cite{IK},
 Kato \cite{Kato}, Simon \cite{Sim}
 serve as a model for us.
 In the paper we use RLS equation by study the stationary scattering problem, when the distance $r{\to}\infty$,  and dynamical  scattering problem, when the time $t{\to}\infty$. We found the connections between these  two problems -ergodic properties. For ordinary differential equations
 the ergodic properties in quantum mechanics were investigated in our papers \cite{Sakh7},\cite{Sakh8}.
  Let us describe the content of the present paper. In section 2 we formulate the results of our previous paper \cite{Sakh4}. Our approach to spectral and scattering problems is essentially based on this previous results. In section 3 we introduce the operator function
  $K(\mu)$ and  investigate properties of this function. Using $K(\mu)$ we constructed the resolvent $(\mathcal{L}-\lambda)^{-1}$ and investigate the spectrum of the operator $\mathcal{L}$. Section 4 is devoted to the construction of the Green's function $G(r,s,\lambda)$ and the  eigenfunction $\phi(r,\lambda)$ that corresponds to the absolutely continuous spectrum. We investigate the properties of $G(r,s,\lambda)$ and  $\phi(r,\lambda)$. In section 5 we investigate the wave operators $W_{\pm}(\mathcal{L},\mathcal{L}_{0})$ and the scattering operator $S(\lambda)$ in the energy representation. We prove that the wave operators $W_{\pm}(\mathcal{L},\mathcal{L}_{0})$ are complete. In section 6 we obtain the following result:\\  The operator
  $S(\lambda)-I$ belongs to the Hilbert-Schmidt class.\\
  It is well-known (see \cite{Gr}) that this result is important in the stationary scattering theory when potential $V(r)$ has the form $V(r)=V(|r|)$. Section 7 is dedicated to ergodic problems in relativistic
  quantum mechanics. It is interesting that our methods are useful for LS equation too. The ergodic problems for LS equation we consider in the last section 8.

  \section{RLS equation in the integral form}
1. Let us write the Dirac equation (see \cite{AB})
\begin{equation}i\frac{\partial}{\partial{t}}u(r,t)=\mathcal{L}u(r,t),\label{2.1}\end{equation}
where $u(r,t)$ is $4\times1$ vector function, $r=(r_1,r_2,r_3)$. The operators $\mathcal{L}$ and $\mathcal{L}_{0}$ are defined by the relations
\begin{equation}\mathcal{L}u=[-e\nu(r)I_{4}+m\beta+{\alpha}({p}+e{A}(r))]u,\quad
\mathcal{L}_{0}u=(m\beta+{\alpha}{p})u. \label{2.2}\end{equation}
Here $p=-igrad$, $\nu$ is a scalar potential, $A$ is a vector potential, $(-e)$ is the electron charge.
 Now let us define $\alpha=[\alpha_1,\alpha_2,\alpha_3]$. The matrices $\alpha_k$  are the $4{\times}4$ matrices of the forms
\begin{equation}\alpha_s=\left(
                           \begin{array}{cc}
                             0 & \sigma_s \\
                             \sigma_s & 0 \\
                           \end{array}
                         \right), \quad s=1,2,3,\label{2.3}\end{equation}
where
\begin{equation}\sigma_1=\left(
                           \begin{array}{cc}
                             0 & 1 \\
                             1 & 0 \\
                           \end{array}
                         \right),\quad
\sigma_2=\left(
                           \begin{array}{cc}
                             0 & -i \\
                             i & 0 \\
                           \end{array}
                         \right),\quad
\sigma_3=\left(
                           \begin{array}{cc}
                             1 & 0 \\
                             0 & -1 \\
                           \end{array}
                         \right).\label{2.4}\end{equation}
The matrices $\beta$ and $I_{2}$ are defined by the relations
\begin{equation} \beta=\left(
                         \begin{array}{cc}
                            I_2 & 0 \\
                           0 & -I_2 \\
                         \end{array}
                       \right),\quad
I_2= \left(
                         \begin{array}{cc}
                          1 & 0 \\
                           0 & 1 \\
                         \end{array}
                       \right).\label{2.5}\end{equation}
2. We consider separately the unperturbed Dirac equation \eqref{2.1}, \eqref{2.2}, when $\nu(r)=0$ and $A(r)=0$. The inverse Fourier transform is defined by the relation
\begin{equation}
\Phi(q)=F^{-1} u(r)=(2\pi)^{-3/2}\int_{R^3}e^{-iqr}u(r)dr. \label{2.6}\end{equation}The Fourier transform has the form
\begin{equation} u(r)=F\Phi(q)=(2\pi)^{-3/2}\int_{Q^3}e^{iqr}\Phi(q)dq.\label{2.7}\end{equation}
In the momentum space the unperturbed  Dirac equation takes the form :
 \begin{equation}i\frac{\partial}{\partial{t}}\Phi(q,t)=H_{0}(q)\Phi(q,t),\quad q=(q_1,q_2,q_3),\label{2.8}\end{equation}
 where
\begin{equation}\mathcal{L}_{0}=F[H_{0}(q)F^{-1}], \label{2.9}\end{equation}
and $H_{0}(q)$ and $\Phi(q,t)$ are matrix functions of order $4{\times}4$ and  $4{\times}1$
respectively.  The matrix $H_{0}(q)$ is defined by the relation
\begin{equation} H_{0}(q)=\left[ \begin{array}{cccc}
                          m & 0 & q_3 & q_1-iq_2 \\
                          0 & m & q_1+iq_2  & -q_3 \\
                          q_3 & q_1-iq_2  & -m & 0 \\
                          q_1+iq_2  & -q_3 & 0 & -m
                        \end{array}\right] .
\label{2.10}\end{equation}
The eigenvalues $\lambda_{k}$ and the corresponding eigenvectors
$g_k$ of  $H_{0}(q)$ are important in our theory.  We have (see \cite{Sakh4}):
\begin{equation}
\lambda_{1,2}=-\sqrt{m^2+|q|^{2}},\quad \lambda_{3,4}=\sqrt{m^2+|q|^{2}},
\label{2.11}\end{equation}where  $|q|:=\sqrt{q_{1}^2+q_{2}^2+q_{3}^2}$.
\begin{equation}
g_1=\begin{bmatrix}(-q_1+iq_2)/(m+\lambda_3) \\ q_3/(m+\lambda_3) \\ 0 \\ 1\end{bmatrix}, \quad
g_2=\begin{bmatrix}-q_3/(m+\lambda_3) \\ (-q_1- iq_2)/(m+\lambda_3) \\1 \\ 0 \end{bmatrix},
\label{2.12}
\end{equation}
\begin{equation}g_3=\begin{bmatrix} (-q_1+iq_2)/(m-\lambda_3)\\ q_3/(m-\lambda_3)\\ 0 \\ 1\end{bmatrix},
\quad
g_4= \begin{bmatrix} -q_3/(m-\lambda_3) \\ (-q_1-iq_2)/(m-\lambda_3) \\ 1 \\ 0\end{bmatrix}.
\label{2.13}\end{equation}
It follows from \eqref{2.10} and \eqref{2.11} that
\begin{equation}
(H_{0}(q)-\lambda)^{-1}=H_{0}^{-1}(q)+H_{0}^{-1}(q)\frac{\lambda^2}
{\lambda_{1}^{2}(|q|)-\lambda^{2}}+\frac{\lambda}{\lambda_{1}^{2}(|q|)-\lambda^{2}}.\label{2.14}\end{equation}
 This equality is valid at all
$\lambda{\notin}E,$ where $E=(-\infty,-m]{\bigcup}[m,+\infty)$.\\
3. \emph{Now we will introduce a relativistic analogue of the Lippmann-Schwinger equation
(RLS integral equation) which was constructed in the paper \cite{Sakh4}}.\\
To do it we consider the expression
\begin{equation}B_{+}(r,\lambda)=Q(r)+(2\pi)^{3/2}\lambda^{2}Q(r)\ast{J_{+}(r,\lambda)}+
{\lambda}J_{+}(r,\lambda),\label{2.15}\end{equation}
where  $\lambda=\overline{\lambda},\, |\lambda|>m$ and
\begin{equation}J_{+}(r,\lambda)=\sqrt{\frac{\pi}{2}}e^{{}i\varkappa|r|}/|r|,\quad \lambda>m,\label{2.16}\end{equation}
\begin{equation}J_{+}(r,\lambda)=\sqrt{\frac{\pi}{2}}e^{{-}i\varkappa|r|}/|r|,\quad \lambda<-m.\label{2.17}\end{equation}
The matrix function Q(r) has the form
\begin{equation} Q(r)=F[H_{0}^{-1}(q)]=\sqrt{\pi/2}e^{-m|r|}[m\beta+i(m+1/|r|)r\alpha/|r|]/|r|.\label{2.18}\end{equation}
Here $r\alpha=r_1\alpha_1+r_2\alpha_2+r_3\alpha_3$, matrices $\alpha_k$ and $\beta$ are defined by the relations  \eqref{2.3}-\eqref{2.5}.
We note that
\begin{equation} F(r)\ast{G(r)}=\int_{R^{3}}F(r-v)G(v)dv \label{2.19}\end{equation} is the convolution of $F(r)$ and $G(r)$.
Now we can write the equation \cite{Sakh4}:
\begin{equation}\phi(r,k,n)=e^{ik{\cdot}r}\hat{g}_{n}(k)-(2\pi)^{3/2}\int_{R^{3}}B_{+}(r-s,\lambda)V(s)\phi(s,k,n)ds,
\label{2.20}\end{equation} where
\begin{equation}V(r)=-e{\nu}(r)I_4+e{\alpha}A(r).\label{2.21}\end{equation}
Here the vectors $g_{n}(k)$ are defined by the relations \eqref{2.12}, \eqref{2.13} and
$\hat{g}_{n}(k)=g_{n}(k)/\|g_{n}\|,\, k{\cdot}r$ is the scalar product of $k$ and $r$.   \\
\begin{Rk}\label{Remark 2.1}Taking into account \eqref{2.11} we have $\lambda=-\sqrt{k^2+m^2},$ when
$n=1,2$, and   $\lambda=\sqrt{k^2+m^2},$ when
$n=3,4$.\end{Rk}
\emph{ Equation \eqref{2.20} (RLS equation) is relativistic analogue of the Lippmann-Schwinger equation.}\\
4. Further we assume that the matrix $V(r)$ is self-adjoint,
\begin{equation}V(r)=V^{\star}(r).\label{2.22}\end{equation}
 Hence $V(r)$ can be represented in the form
\begin{equation}V(r)=U(r)D(r)U^{\star}(r),\label{2.23}\end{equation}
where $U(r)$ is an unitary matrix, $D(r)$ is a diagonal matrix
\begin{equation}D(r)=diag(d_{1}(r),d_{2}(r),d_{3}(r),d_{4}(r)).\label{2.24}\end{equation}
Let us introduce the diagonal matrices
\begin{equation}D_{1}(r)=diag(|d_{1}(r)|^{1/2},|d_{2}(r)|^{1/2},|d_{3}(r)|^{1/2},|d_{4}(r)|^{1/2})
\label{2.25}\end{equation}
and
\begin{equation}W(r)=diag(sgnd_{1}(r),sgnd_{2}(r),sgnd_{3}(r),sgnd_{4}(r)).\label{2.26}\end{equation}
Formulas \eqref{2.23}-\eqref{2.26} imply that
\begin{equation}V(r)=V_{1}(r)W_{1}(r)V_{1}(r),\label{2.27}\end{equation}
where
\begin{equation}V_{1}(r)=U(r)D_{1}(r)U^{\star}(r),\quad W_{1}(r)=U(r)W(r)U^{\star}(r)\label{2.28}\end{equation}
It is easy to see that
\begin{equation}\|V_{1}(r)\|^{2}=\|V(r)\|,\quad \|W_{1}(r)\|=1.\label{2.29}\end{equation}
 \emph{Modified RLS integral equation.}\\
If the $4{\times}1$ matrix function $\phi(r,k,n)$ is a solution of RLS equation, then the matrix function $\psi(r,k,n)=V_{1}(r)\phi(r,k,n)$
is a solution of following  modified RLS integral equation:
\begin{equation}\psi(r,k,n)=e^{ik{\cdot}r}V_{1}(r)\hat{g}_{n}(k)-(2\pi)^{3/2}B_{+}(\lambda)\psi(r,k,n),
\label{2.30}\end{equation}
where
\begin{equation}B_{+}(\lambda)f=\int_{R^{3}}V_{1}(r)B_{+}(r-s,\lambda)V_{1}(s)W_{1}(s)f(s)ds.
\label{2.31}\end{equation} We note that the operators $B_{+}(\lambda)$ act   in the Hilbert space
$L^{2}(R^{3})$. (We say that the matrix belongs to the Hilbert space
$L^{2}(R^{3})$ if the every element of the matrix belongs to the Hilbert space
$L^{2}(R^{3})$.)
We have proved the following result \cite{Sakh4}.
\begin{Tm}\label{Theorem 2.2}If
the function  $\|V(r)\|$ is bounded and belongs to the space $L^{1}(R^{3})$
then the operator $B_{+}(\lambda)$ is compact.
\end{Tm}
5. Let us introduce the following definition \cite{Sakh4}
\begin{Dn}\label{Definition 2.3}We say that $\lambda{\in}E$ is an exceptional value if the equation $[I+(2\pi)^{3/2}B_{+}(\lambda)]\psi=0$ has nontrivial solution in the space $L^2(R^3)$.\end{Dn}
 We denote by $\mathcal{E}_{+}$ the set of exceptional points and we denote by $E_{+}$ the set
 of such points $\lambda$ that $\lambda{\in}E,\quad \lambda{\notin}\mathcal{E}_{+}.$\\
 We have \cite{Sakh4}:
\begin{La}\label{Lemma 2.4}Let conditions of Theorem 2.2 be fulfilled. If $\lambda{\in}E_{+}$, then equation \eqref{2.30} has one and only one
solution $\psi(r,k,n)$ in $L^2(R^3)$. \end{La}
\begin{Cy}\label{Corollary 2.5}Let conditions of Theorem 2.2 be fulfilled. If $\lambda{\in}E_{+}$
, then equation \eqref{2.20} has one and only one
solution $\phi(r,k,n)$ which satisfies the condition $V_{1}(r)\phi(r,k,n){\in}L^2(R^3)$. \end{Cy}
5. The following Theorem gives the connection between spectral and scattering results \cite{Sakh4}.
\begin{Tm}\label{Theorem 2.6}Let the $V(r)=V^*(r)$ and the function $\|V(r)\|$ be bounded and belong to the space $L^{1}(R^3)$. If
\begin{equation}\|V_1(r)\|=O(|r|^{-3/2}),\quad |r|{\to}\infty,\label{2.32}\end{equation}
then the solution $\phi(r,k,n,\lambda)$ of RLS equation \eqref{2.20} has the form
\begin{equation}\phi(r,k,n,\lambda)=e^{ik{\cdot}r}\hat{g}_{n}(k)+\frac{e^{i\varkappa|r|}}{|r|}
f(\omega,\omega^{\prime},n,\lambda)
+o(1/|r|), \label{2.33}\end{equation}
where  $|r|{\to}\infty$ , $\lambda{\in}E_{+}, \,\omega=r/|r|,\, \omega^{\prime}=k/|k|$ and
\begin{equation}f(\omega,\omega^{\prime},n,\lambda)=-2\pi^2{\lambda}\int_{R^3}e^{-i{\varkappa}s{\cdot}\omega}
V(s)\phi(s,k,n,\lambda)ds.  \label{2.34}\end{equation}
 \end{Tm}
(We note that $k$ is a vector from the space $R^3$  and $\varkappa=|k|=\sqrt{\lambda^2-m^2}.$)
 \begin{Dn}\label{Definition 2.7} The $4{\times}1$ vector functions $f(\omega,\omega^{\prime},n,\lambda)$ are named the relativistic scattering
amplitudes.\end{Dn}

6. Now we shall formulate the connection between solutions of the equation
\begin{equation}\mathcal{L}\phi=\lambda\phi\label{2.35}\end{equation}
and the solutions of the  RLS  equation \eqref{2.20}(see \cite{Sakh4}).
\begin{Tm}\label{Theorem 2.8}Let the conditions of Theorem 2.2 be fulfilled and let  the function $\phi(r,k,n)$ be a solution of RLS equation \eqref{2.20} such that \\ $V_{1}(r)\phi(r,k,n){\in}L^{2}((R^3).$ If $\lambda{\in}E_{+}$ then the vector function $\phi(r,k.n)$
is the solution of the equation \eqref{2.35}
in the distributive sense.
\end{Tm}
\begin{Cy}\label{Corollary 2.9} Let potential $V(r)$ satisfy the conditions of Theorem 2.2. If
$\lambda{\in}E$ is an eigenvalue of the corresponding operator $\mathcal{L}$, then $\lambda{\in}\mathcal{E}_{+}$.\end{Cy}
\section{Definition and properties of the holomorphic operator function $K(\mu)$.}
1.\emph{ Let us consider the case $V(r)=0$ separately.}\\
The function $B_{+}(r,\lambda)$ can be determined not only under the condition $\lambda{\in}E.$
Taking into account \eqref{2.15} we have
\begin{equation}
  B_{+}(r,\mu)= Q(r)+(2\pi)^{3/2}\mu^{2}Q(r)\ast{J_{+}(r,\mu)}+
{\mu}J_{+}(r,\mu),\label{3.1}\end{equation}
where $\Im{\mu}>0,\, \varkappa=\sqrt{\mu^2-m^2},\, \Im{\varkappa}>0$ and
\begin{equation}J_{+}(r,\mu)=\sqrt{\frac{\pi}{2}}e^{i\varkappa|r|}/|r|.
\label{3.2}\end{equation}
We introduce the operator
\begin{equation}{B}_{0}(\mu)f=\int_{R^3}B_{+}(r-s,\mu)f(s)ds,\label{3.3} \end{equation}
where $f(r)$ is  $4{\times}1$ vector function with elements in the space $L^2(R^3)$. It is easy to see that\\
\begin{equation}\int_{R^3}\|B_{+}(r,\mu)\|dr<\infty,\quad \Im{\varkappa}>0.\label{3.4}\end{equation}
Hence
the operator ${B}_{0}(\mu)$
is bounded in the space $L^2(R^3)$ (see \cite{Sakh6}, p.19).
It follows from \eqref{2.14} that
\begin{equation}(\mathcal{L}_{0}-\mu)^{-1}f=(2\pi)^{3/2}\int_{R^3}B_{+}(r-s,\mu)f(s)ds,\quad \Im{\mu}>0,\quad \Im{\varkappa}>0 .\label{3.5}\end{equation}
2. Now we consider the operators :
\begin{equation}B_{+}(\mu)f=\int_{R^{3}}V_{1}(r)B_{+}(r-s,\mu)V_{1}(s)W_{1}(s)f(s)ds,\quad
\Im{\mu}>0,\quad \Im{\varkappa}>0
\label{3.6}\end{equation}
and
\begin{equation}K(\mu)=I+(2\pi)^{3/2}B_{+}(\mu).\label{3.7}\end{equation}
\begin{La}\label{Lemma 3.1}Let the conditions of theorem 2.2 be fulfilled. Then the operator function $K(\mu)$ has the following properties in the region $\mu{\notin}E:$\\
1)The operator function $K(\mu)-I$ is holomorphic.\\
2) The operator $K(\mu)-I$ is compact.\\
3) For some $\mu_{0}$ the operator $K(\mu_{0})$ has a bounded inverse operator.\end{La}
\emph{Proof.} The assertion 1) follows from \eqref{3.1} and \eqref{3.5}.  The assertion 2) is proved in the paper (\cite{Sakh4}, Theorem 3.1) for the case when $\lambda{\in}E$. The proof is correct for the case when $\lambda{\notin}E,\,\Im{\varkappa}>0$. To prove the assertion 3) we use the the following property of the self-adjoint operators
\begin{equation}\|(\mathcal{L}_{0}-\mu)^{-1}\|{\leq}|\Im{\mu}|^{-1}.\label{3.8}
\end{equation}Then according to \eqref{3.5} - \eqref{3.8} we obtain
\begin{equation}\|B_{+}(\mu)\|{\to}0,\quad K(\mu){\to}I, \quad \Im{\mu}{\to}\infty.
\label{3.9}\end{equation}
 The assertion 3) follows directly from \eqref{3.9}. The lemma is proved.\\
 It follows from Lemma 3.1, that all conditions of Gohberg theorem (see \cite{Br}, Appendix II)
 are fulfilled. Then we have:
 \begin{Pn}\label{Proposition 3.2} Let conditions of Theorem 2.2 be fulfilled. If $\mu_{0}{\notin}E$, then either
 the operator function $K^{-1}(\mu)$ is holomorphic in the point $\mu_{0}$ or the point $\mu_{0}$ is a pole of the operator function $K^{-1}(\mu)$.\end{Pn}
 We need the following assertion.
\begin{Tm}\label{Theorem 3.3}Let the conditions of  Theorem 2.2 be fulfilled. Then
\begin{equation}(\mathcal{L}-\mu)^{-1}-(\mathcal{L}_{0}-\mu)^{-1}=-(\mathcal{L}_{0}-\mu)^{-1}V_1W_1
K(\mu)^{-1}V_1(\mathcal{L}_{0}-\mu)^{-1}, \label{3.10}\end{equation}
where $\mu{\notin}E$.\end{Tm}
\emph{Proof.} Let us write the relation
\begin{equation}(\mathcal{L}-\mu)f=(\mathcal{L}_{0}-\mu+V(r)))f=g.\label{3.11}\end{equation}
It follows from \eqref{3.11} that
\begin{equation}(\mathcal{L}-\mu)^{-1}=K_{1}(\mu)^{-1}(\mathcal{L}_{0}-\mu)^{-1},\label{3.12}\end{equation}
where
\begin{equation}K_{1}(\mu)=I+(\mathcal{L}_{0}-\mu)^{-1}V.\label{3.13}\end{equation}
It follows from \eqref{3.12} and \eqref{3.13} that
\begin{equation}(\mathcal{L}-\mu)^{-1}- (\mathcal{L}_{0}-\mu)^{-1}=-(\mathcal{L}_{0}-\mu)^{-1}VK_{1}^{-1}(\mu)(\mathcal{L}_{0}-\mu)^{-1}.
\label{3.14}\end{equation}
Taken into account equality $V=V_1W_1V_1$ and relations \eqref{3.6}, \eqref{3.7} and \eqref{3.13} we have
\begin{equation}V_{1}K_{1}(\mu)=K(\mu)V_1,\quad V_{1}K_{1}^{-1}(\mu)=K^{-1}(\mu)V_1.\label{3.15}\end{equation}
The assertion of the Theorem follows from \eqref{3.14} and \eqref{3.15}.
\begin{Rk}\label{Remark 3.4} The equality of type \eqref{3.10} was proved in the book (\cite{Sim}, p.73).
\end{Rk}
The following assertion is true (see \eqref{2.9} and \eqref{2.10}):
\begin{Pn}\label{Proposition 3.5} The operator $\mathcal{L}_{0}$ has only absolutely continuous spectrum,
which coincides with the set E.
\end{Pn}
\begin{Cy}\label{Corollary 3.6}Let the conditions of theorem 2.2 be fulfilled.
Then\\
1) The operator $\mathcal{L}$ has only discrete spectrum in the interval
$-m<\lambda<m$.\\
2) The discrete spectrum $-m<\lambda_{n}<m$ of the  operator $\mathcal{L}$ has no limit points in the interval $-m<\lambda<m$.\\
3)The set $\mathcal{E}_{+}$ is   closed and has  Lebesgue  measure equal to zero. \end{Cy}
\emph{Proof.} Assertions 1)and 2) follows from relation \eqref{3.10}, Proposition 3.2 and Proposition 3.4.\\
Assertion 3) follows from Lemma 3.1 in the same way as in the case  of
 classical
Lippmann-Schwinger equation (\cite{RS}, p.115).\\

\section{Green matrix-function and eigenfunction expansion}

1. In this section we construct the Green matrix function of the operator $\mathcal{L}$ (see \eqref{2.2}). Multiplying \eqref{3.10} on the left $V_{1}$  or on the right $V_{1}W_{1}$ and using Theorem 2.2 we conclude.
\begin{Cy}\label{Corollary 4.1}  Let conditions of Theorem 2.2 be fulfilled and $\Im{\mu}>0$. Then
the operators  $V_{1}(\mathcal{L}-\mu)^{-1}$  and $(\mathcal{L}-\mu)^{-1}V_{1}$ are compact.\end{Cy}
Let us consider the integral equation
\begin{equation}G(r,s,\mu)=G_{0}(r,s,\mu)-\int_{R^3}G_{0}(r,z,\mu)V(z)G(z,s,\mu)dz.\label{4.1}
\end{equation}
In view of \eqref{3.1} and \eqref{3.5} the relation
\begin{equation}G_{0}(r,s,\mu)=(2\pi)^{3/2}B_{+}(r-s,\mu),\,\Im{\mu}>0 \label{4.2}\end{equation} holds.
\begin{Tm}\label{Theorem 4.2.} Let conditions of Theorem 2.2 be fulfilled and $\Im{\mu}>0,\quad \Im{\varkappa>0}$.
Then\\
1)Equation \eqref{4.1} has one and only one solution such that $V_{1}(r)G(r,s,\mu){\in}L^{2}(R^3)$,\\where $s$ and $\mu$ are fixed.\\
2) The solution $G(r,s,\mu)$ of equation \eqref{4.1} is the  Green matrix function of the operator
$\mathcal{L}$, i.e. the following equality
\begin{equation}R_{\mu}f=(\mathcal{L}-\mu)^{-1}f(r)=\int_{R^3}G(r,s,\mu)f(s)ds.\label{4.3}\end{equation}
is valid.\\
3) The Green   matrix function $G(r,s,\mu)$ belongs to the space $L^{1}(R^3)$, where $s$ and $\mu$ are fixed.\\
4) The Green   matrix function $G(r,s,\mu)$ belongs to the space $L^{1}(R^3)$, where $r$ and $\mu$ are fixed.\\
\end{Tm}
\emph{Proof.} The equality \eqref{4.1} is equivalent to the relation
\begin{equation}(\mathcal{L}-\mu)^{-1}-(\mathcal{L}_{0}-\mu)^{-1}=
(\mathcal{L}_{0}-\mu)^{-1}V(\mathcal{L}-\mu)^{-1}.\label{4.4}\end{equation}
 Then the equality \eqref{4.3}  holds, i.e. the assertion 1) and 2) are proved.
 To prove 3) we need to show
 \begin{equation}\int_{R^3}\int_{R^3}\|G_{0}(r,z,\mu)V(z)G(z,s,\mu)\|dzdr<\infty.\label{4.5}
 \end{equation} We note that
$\int_{R_3}\|G_{0}(r,z,\mu)\|dr$  is finite independent of z (see \eqref{4.2}). Using relations
\begin{equation}\|V_{1}(z)\|{\in}L^{2}(R^3), \quad \|V_{1}(z)G(z,s,\lambda)\|{\in}L^{2}(R^3)
\label{4.6}\end{equation}
we obtain
 \begin{equation}\int_{R^3}\|V(z)G(z,s,\mu)\|dz<\infty.\label{4.7}
 \end{equation}
Now assertion 3)  follows directly from \eqref{4.1} and inequality
\begin{equation}\int_{R^3}\|B_{+}(r,\mu)\|dr<\infty.\nonumber
 \end{equation}
The assertion 4) follows from assertion 3) and equality $G^{\star}(r,s,\mu)=G(s,r,\overline{\mu}).$\\ The Theorem is proved.\\
2. Taking into account Theorem 4.2 we introduce the $4{\times}4$ matrix function
\begin{equation} g(r,k,\mu)=(2\pi)^{-3/2}\int_{R^3}G(r,s,\mu)e^{is{\cdot}k}ds=
F (G(r,s,\mu)).\label{4.8}\end{equation}
In particular, we have
\begin{equation} g_{0}(r,k,\mu)=(2\pi)^{-3/2}\int_{R^3}G_{0}(r,s,\mu)e^{is{\cdot}k}ds=F(G_{0}(r,s,\mu)).\label{4.9}
\end{equation}
It follows from \eqref{4.2} that
\begin{equation}F^{-1}(G_{0}(r.s,\mu))=e^{ir{\cdot}k}(2\pi)^{3/2}F^{-1}(B_{+}(u,\mu)),\quad r-s=u.\label{4.10}\end{equation}
Taking into account the relation (see \cite{Sakh4}, formula \eqref{2.17}):
\begin{equation}B_{+}(r,\lambda)=F[H_{0}(q)-(\lambda+i0)]^{-1},\nonumber\end{equation}
 we have
\begin{equation}B_{+}(r,\mu)=F[H_{0}(q)-\mu]^{-1},\quad \Im{\mu}>0.\label{4.11}\end{equation}
It follows from \eqref{4.9}-\eqref{4.11} that
\begin{equation} g_{0}(r,k,\mu)=(2\pi)^{3/2}e^{ik{\cdot}r}[H_{0}(k)-\mu]^{-1}.\label{4.12}\end{equation}
We apply to \eqref{4.1} the Fourier transform $F$ with respect to $s$:
\begin{equation}h(r,k,\mu)=e^{ik{\cdot}r}-\int_{R^3}G_{0}(r,z,\mu)V(z)h(z,k,\mu)dz,\label{4.13}
\end{equation}
where
\begin{equation}h(r,k,\mu)=(2\pi)^{-3/2}g(r,k,\mu)[H_{0}(k)-\mu].\label{4.14}\end{equation}
So, we have
\begin{equation}h(r,k,\mu)=e^{ik{\cdot}r}-\int_{R^3}G_{0}(r,z,\mu)W_{1}(z)V_{1}(z)p(z,k,\mu)dz,\label{4.15}\end{equation}
where $p(r,k,\mu)=V_{1}(r)h(r,k,\mu)$ is unique solution of equation
\begin{equation}p(r,k,\mu)=V_{1}(r)e^{ikr}-\int_{R^3}V_{1}(r)G_{0}(r,z,\mu)W_{1}(z)V_{1}(z)p(z,k,\mu)dz,\label{4.16}\end{equation}
Let us consider the vector functions
\begin{equation}\hat{h}_{n}(r,k,\mu)=h(r,k,\mu)\hat{g}_{n}(k), \quad
\hat{p}_{n}(r,k,\mu)=p(r,k,\mu)\hat{g}_{n}(k)\label{4.17}\end{equation}
Then equalities \eqref{4.15} and \eqref{4.16} can be rewritten in the forms:
\begin{equation}\hat{h}_{n}(r,k,\mu)=e^{ik{\cdot}r}\hat{g}_{n}(k)-
\int_{R^3}G_{0}(r,z,\mu)W_{1}(z)V_{1}(z)\hat{p}_{n}(z,k,\mu)dz\label{4.18}\end{equation}
and
\begin{equation}\hat{p}_{n}(r,k,\mu)=e^{ik{\cdot}r}V_{1}(r)\hat{g}_{n}(k)-(2\pi)^{3/2}B_{+}(\mu)\hat{p}_{n}(r,k,\mu),
\label{4.19}\end{equation}
where
\begin{equation}B_{+}(\mu)f=\int_{R^{3}}V_{1}(r)B_{+}(r-s,\mu)V_{1}(s)W_{1}(s)f(s)ds.
\label{4.20}\end{equation}
3. Now  we use the results of point  2 and study  the absolutely continuous spectrum of the
operator $\mathcal{L}$.
 We assume that $f(r){\in}C_{0}^{\infty}$ and consider the integrals
\begin{equation}\Phi_{n}(k,\mu)=(2\pi)^{-3/2}\int_{R^3}\hat{h}^{*}_{n}(r,k,\mu)f(r)dr, \label{4.21}\end{equation}
\begin{equation}\widetilde{f}_{n}(k)=(2\pi)^{-3/2}
\int_{R^3}\phi^{*}_{n}(r,k)f(r)dr.
\label{4.22}\end{equation}
\begin{La}\label{Lemma 4.3}  Let conditions of Theorem 2.2 be fulfilled and $\Im{\mu}>0$. \\
If $0<m<\alpha<\beta$ and $[\alpha,\beta]{\in}E_{+}$,
 then  $\Phi_{n}(k,\mu),(\Re{\mu}=\varkappa,\quad |k|=\varkappa,\quad n=3,4)$  can be extended uniformly continuous to the region:\\
$$ \sqrt{\varkappa^2+m^2}{\in}[\alpha,\beta].$$\\
and  the following equality
\begin{equation} \Phi_{n}(k,\varkappa)=\widetilde{f}_{n}(k),\quad (|k|=\varkappa),\quad n=3,4\label{4.23}\end{equation}
is valid.
\end{La}
\emph{Proof.} If $k{\in}R^3$ and $|k|=\varkappa$, then the equation \eqref{4.19} is identical to the equation \eqref{2.30}. Hence, in this case we have: $\hat{p}_{n}(r,k,\varkappa)=\psi_{n}(r,k)$    and $\hat{h}_{n}(r,k,\varkappa)=\phi_{n}(r,k)$. The assertion of the lemma follows directly
from the  relations \eqref{4.21} and \eqref{4.22}.\\
In the same way as Lemma 4.3 we obtain
\begin{La}\label{Lemma 4.4} Let conditions of Theorem 2.2 be fulfilled and $\Im{\mu}>0$.\\
 If $ 0>-m>\beta>\alpha$ and $[\alpha,\beta]{\in}E_{+},$
 then  $\Phi_{n}(k,\mu),(\Re{\mu}=-\varkappa<0,\quad|k|=\varkappa,\quad n=1,2)$  can be extended uniformly continuous to the region\\
$$|k|=\varkappa, \quad -\sqrt{\varkappa^2+m^2}{\in}[\alpha,\beta]$$\\
and the following equality
\begin{equation} \Phi_{n}(k,-\varkappa)=\widetilde{f}_{n}(k),\quad (|k|=\varkappa,\quad n=1,2)\label{4.24}\end{equation}
is valid.
\end{La}
\begin{Dn}\label{Definition 4.5} We denote by $\mathcal{H}_{ac}$ the absolutely continuous invariant subspace with respect to $\mathcal{L}$ and by $P_{ac}$ orthogonal
projector on  $\mathcal{H}_{ac}$.\end{Dn}
\begin{Dn}\label{Definition 4.6} We denote by $\mathcal{H}_{\alpha}$ the maximal invariant subspace with respect to $\mathcal{L}$ such that the spectrum of $\mathcal{L}$ belongs to the set $(-\infty,\alpha]$. We denote  by $P_{\alpha}$ orthogonal
projector on  $\mathcal{H}_{\alpha}$.\end{Dn}
Consider the $4{\times}4$ matrix
\begin{equation} Z_{0}(k)=[\hat{g}_{1}(k),\hat{g}_{2}(k),\hat{g}_{3}(k),\hat{g}_{4}(k)],
\label{4.25}\end{equation}
\begin{Pn}\label{Proposition 4.7}The matrix $Z_{0}(k)$ is unitary and
\begin{equation}H_{0}(k)=Z_{0}(k)\mathcal{D}(k)Z_{0}^{\star}(k),\label{4.26}\end{equation}
and the diagonal matrix $\mathcal{D}(k)$ is defined by the relation
\begin{equation}\mathcal{D}(k)=Diag[\lambda_1(k),\lambda_2(k),\lambda_3(k),\lambda_4(k)].
\label{4.27}\end{equation}
\end{Pn}
\begin{La}\label{Lemma 4.8} Let conditions of Theorem 2.2 be fulfilled. If \\ $0<m<\alpha<\beta,\,[\alpha,\beta]{\in}E_{+},$
then
\begin{equation}\|(P_{\beta}-P_{\alpha})f\|^{2}=\int_{a<\varkappa<b}[|\widetilde{f}_{3}(k)|^{2}+
|\widetilde{f}_{4}(k)|^{2}]dk,
\label{4.28}
\end{equation}
where $a=\sqrt{\alpha^2-m^2},\quad b=\sqrt{\beta^2-m^2},\quad \varkappa=|k|.$
\end{La}
\emph{Proof.}
Parseval equality for Fourier transform \eqref{4.8} implies that
\begin{equation}\int_{R^3}G(s,z,\mu)G^{*}(r,z,\mu)dz=
\int_{R^3}g(s,k,\mu)g^{*}(r,k,\mu)dk,\label{4.29}\end{equation}
where $\mu=\lambda+ i\varepsilon,\, \lambda=\overline{\lambda},\,\varepsilon>0.$
We introduce  the notation
\begin{equation}R_{0}(\mu)=[D(k)-\mu]^{-1}.
\label{4.30}\end{equation}
Relations \eqref{4.14}, \eqref{4.26} and \eqref{4.30} imply  that
\begin{equation}\int_{R^3}G(s,z,\mu)G^{*}(r,z,\mu)dz=
\int_{R^3}\hat{H}(s,k,\mu)T(\mu,\overline{\mu})\hat{H}^{*}(r,k,\mu)dk.
\label{4.31}\end{equation}
where
\begin{equation}\hat{H}(s,k,\mu)=h(s,k,\mu)Z_{0}(k),\quad T(\mu,\overline{\mu})=(2\pi)^{3}R_{0}(\mu)R_{0}(\overline{\mu}).\label{4.32}\end{equation}
We multiply both sides \eqref{4.29} of left hand side by $f^{\star}(s)$ and
of right hand side by $f(r)$ and integrate by s and by r:
\begin{equation} \int_{R^3}\Psi^{*}(z,\mu)\Psi(z,\mu)dz=
\int_{R^3}\Phi^{*}(k,\mu)T(\mu,\overline{\mu})\Phi(k,\mu)dk,\label{4.33}
\end{equation}where
\begin{equation}\Psi(z,\mu)=\int_{R^3}G^{*}(r,z,\mu)f(r)dr,\quad
\Phi(k,\mu)=\int_{R^3}\hat{H}^{*}(r,k,\mu)f(r)dr.\label{4.34}\end{equation}
Using relation \eqref{4.3} and equality
$G^{\star}(r,s,\mu)=G(s,r,\mu)$
we have
\begin{equation}R_{\overline{\mu}}f=\Psi(r,\mu).\label{4.35}\end{equation}
It follows from \eqref{4.33}   and \eqref{4.35} that
\begin{equation} (R_{\overline{\mu}}f,R_{\overline{\mu}}f)=
\int_{R^3}\Phi^{*}(k,\mu)T(\mu,\overline{\mu})\Phi(k,\mu)dk,\label{4.36}
\end{equation}
Relation \eqref{4.36} implies that
\begin{equation} (f,[R_{\mu}-R_{\overline{\mu}}]f)=
2i\varepsilon[
\int_{R^3}\Phi^{*}(k,\mu)T(\mu,\overline{\mu})\Phi(k,\mu)dk].\label{4.37}\end{equation}
Now we use the fundamental relation (see \cite{St}, p.183).
\begin{equation}(f,[(P_{\beta}+P_{\beta-0})-(P_{\alpha}+P_{\alpha-0})]f)=
\frac{1}{i\pi}{lim}\int_{\alpha}^{\beta}(f,[R_{\lambda+i\varepsilon}-R_{\lambda-i\varepsilon}]f)d\lambda,
\label{4.38}\end{equation} where $\varepsilon{\to}+0$.
Further we need the well-known relation (see \cite{Tit}, p.31)
\begin{equation}J_{(\alpha,\beta)}=\frac{1}{\pi}\lim_{\varepsilon{\downarrow}0}\int_{\alpha}^{\beta}
\frac{\varepsilon}{(c-\lambda)^2+\varepsilon^2}f(\lambda,\varepsilon)d\lambda=f(c,o),\quad \alpha<c<\beta. \label{4.39}\end{equation}
Here $f(\lambda,\varepsilon)$ is a continuous function of $\lambda$ and $\varepsilon$ for $\lambda{\in}[\alpha,\beta]$ and $0{\leq}\varepsilon{\leq}\varepsilon_{0}<\infty.$ It is easy to see that
\begin{equation}J_{(\alpha,\beta)}=0,\quad c{\notin}[\alpha,\beta].\label{4.40}\end{equation}
Let us write one more relation
\begin{equation}\lim_{\varepsilon{\downarrow}0}\int_{\alpha}^{\beta}
\frac{\varepsilon}{\sqrt{(c-\lambda)^2+\varepsilon^2}}f(\lambda,\varepsilon)d\lambda=0. \label{4.41}\end{equation}
Formulas \eqref{4.23}, \eqref{4.37} and \eqref{4.39}-\eqref{4.41} imply
\begin{equation}(f,[(P_{\beta}+P_{\beta-0})-(P_{\alpha}+P_{\alpha-0})]f)=
2\int_{a<\varkappa<b}[|\widetilde{f}_{3}(k)|^{2}+|\widetilde{f}_{4}(k)|^{2}]dk.
\label{4.42}\end{equation}
Letting $\alpha{\to}\beta$ we obtain that
\begin{equation}P_{\beta}=P_{\beta-0}.\label{4.43}\end{equation}
We note that relation \eqref{4.42} is proved for $f{\in}C_{0}^{\infty}$. Hence this relation is
valid for all $f{\in}L^2(R^3)$.The assertion of the lemma follows directly
from \eqref{4.42} and \eqref{4.43}.\\
In the same way as Lemma 4.8 can be proved the following result.
\begin{La}\label{Lemma 4.9} Let conditions of Theorem 2.2 be fulfilled. If  $\alpha<\beta<-m , \,(\alpha,\beta){\in}E_{+},$
then
\begin{equation}\|(P_{\beta}-P_{\alpha})f\|^{2}=\int_{a<\varkappa<b}[|\widetilde{f}_{1}(k)|^{2}+
\widetilde{f}_{2}(k)|^{2}]dk,
\label{4.44}
\end{equation}
and
\begin{equation}P_{\beta}=P_{\beta-0},\label{4.45}\end{equation}
where $b=\sqrt{\alpha^2-m^2},\quad a=\sqrt{\beta^2-m^2},\quad \varkappa=|k|.$
\end{La}
\begin{Tm} \label{Theorem 4.10}Let conditions of Theorem 2.2 be fulfilled. Then\\
1)There exists no eigenvalue of $\mathcal{L}$ which belongs to $E_+$.\\
2)The spectrum of $\mathcal{L}$  on the $E_+$ is absolutely continuous.\\
3) The set of singular spectrum of  $\mathcal{L}$ belongs to $\mathcal{E}_{+}$ and  has  Lebesgue  measure equal to zero. \end{Tm}
\emph{Proof.} The assertion 1) follows from \eqref{4.43} and \eqref{4.45}. To prove assertion 2) we introduce the function
$\sigma_{f}(\beta)=(P_{\beta}f,f)$. Then we have
\begin{equation}\sigma_{f}(\beta)-\sigma_{f}(\alpha)=(f,(P_\beta-P_\alpha)f)=\|(P_\beta-P_\alpha)f\|^{2}
\label{4.46}\end{equation} Relations \eqref{4.28},\eqref{4.44} and \eqref{4.46} imply that the function
$\sigma_{f}(\beta)$ is absolutely continuous when $\beta{\in}E_+.$ Hence the assertion 2) is proved.
It follows from Corollary 3.6 and assertion 2) that singular spectrum of  $\mathcal{L}$ belongs to $\mathcal{E}_{+}$.
We have proved \cite{Sakh4} that  $\mathcal{E}_{+}$ has  Lebesgue  measure equal to zero.
Hence the assertion 3) is valid. The theorem is proved.\\
Let us introduce the domain  $D(R,\varepsilon,k)$  such that
$|k|<R$ and
$|k-k_{0}|>\varepsilon>0$ for all $k_{0}$ satisfying the condition
$\pm\sqrt{k_{0}^{2}+m^2}{\in}\mathcal{E}_{+}$.
 We shall use the following notation:
 \begin{equation} \int_{D}g(k)dk=lim\int_{D(R,\varepsilon,k)}g(k)dk,\quad R{\to}\infty,\quad
 \varepsilon{\to}+0.\label{4.47}\end{equation}
\begin{Pn}\label{Proposition 4.11}Let conditions of Theorem 2.2 be fulfilled. Then
\begin{equation}\|P_{ac}f\|^{2}=\int_{D}[|\widetilde{f}_{1}(k)|^{2}+
|\widetilde{f}_{2}(k)|^{2}+|\widetilde{f}_{3}(k)|^{2}+|\widetilde{f}_{4}(k)|^{2}]dk \label{4.48} \end{equation}
and
\begin{equation}P_{ac}f=(2\pi)^{-3/2}\int_{D}\sum_{p=1}^{4}[\phi_{p}(r,k,)\widetilde{f}_{p}(k)]dk
\label{4.49}\end{equation}
\end{Pn}
\emph{Proof.} The equality \eqref{4.48} follows from \eqref{4.28} and \eqref{4.44}. To prove \eqref{4.49} we consider
 $g(r){\in}C_{0}^{\infty}.$  Using \eqref{4.48} we obtain
\begin{equation}(g,P_{ac}f)=\int_{D}\sum_{p=1}^{4}[\widetilde{g}^{\star}_{p}(k)
\widetilde{f}_{p}(k)]dk
\label{4.50}\end{equation} where (see \eqref{4.22})
\begin{equation}\widetilde{g}_{n}(k)=(2\pi)^{-3/2}
\int_{R^3}\phi^{*}_{n}(r,k)g(r)dr.
\label{4.51}\end{equation}
 Relation \eqref{4.49} follows from \eqref{4.50} and \eqref{4.51}.
Proposition is proved.
\begin{Rk}\label{Remark 4.12} Suppose for simplicity that the corresponding operator
$\mathcal{L}$ has
 no singular and discrete spectra. Then formulas \eqref{4.22} and \eqref{4.44}
 are expansion formulas  in terms of  generalized eigenfunctions $\phi_{n}(r,k)$.\end{Rk}
   \section{Scattering theory,  completeness of wave operators}
1. In this section we use the constructed relativistic Lippman-Schvinger equation
to study the  scattering problems.\\
We introduce the operator function
\begin{equation}\Theta(t)=exp(it\mathcal{L})exp(-it\mathcal{L}_{0}).
\label{5.1}\end{equation}
The wave operators $W_{\pm}(\mathcal{L},\mathcal{L}_0)$ are defined by the relation (see \cite{RS}).
\begin{equation}W_{\pm}(\mathcal{L},\mathcal{L}_0)=\lim_{t{\to}\pm{\infty}}\Theta(t)P_0.\label{5.2}
\end{equation} Here  $P_0$ is orthogonal projector
on the absolutely continuous subspace $G_0$ with respect to the operator $\mathcal{L}_0$.  The
limit in \eqref{5.2} supposed to be  in the sense of strong convergence.
\begin{Rk}\label{Remark 5.1} We consider the case when the operator $\mathcal{L}_{0}$ is defined by the relation \eqref{2.2}. In this case we have: $G_{0}=R^3$ and $P_{0}=I$.\end{Rk}
We have proved the assertion \cite{Sakh4}.
\begin{Tm}\label{Theorem 5.2} If $V(r)=V^{*}(r)$,
the function  $\|V(r)\|$ is bounded, belongs to the space $L^{1}(R^{3})$ and
\begin{equation}\int_{-\infty}^{+\infty}[\int_{|r|>|t|\varepsilon}\|V(r)\|^{2}dr]^{1/2}dt<\infty,\quad \varepsilon>0, \label{5.3}\end{equation}
then the wave operators $W_{\pm}(\mathcal{L},\mathcal{L}_0)$ exist.
\end{Tm}
\begin{Cy}\label{Corollary 5.3}Let the condition  $V(r)=V^{*}(r)$ be fulfilled. If
the function  $\|V(r)\|$ is bounded and
\begin{equation}\|V(r)\|{\leq}\frac{M}{|r|^{\alpha}},\quad |r|{\geq}\delta>0,\quad \alpha>3,
\label{5.4}\end{equation}
then the wave operators $W_{\pm}(\mathcal{L},\mathcal{L}_0)$ exist.
\end{Cy}
\begin{Dn}\label{Definition 5.4}The scattering operator $S(\mathcal{L},\mathcal{L}_0)$ is defined by the relation
\begin{equation}S(\mathcal{L},\mathcal{L}_0)=W_{+}^{*}(\mathcal{L},\mathcal{L}_0)W_{-}(\mathcal{L},\mathcal{L}_0).
\label{5.5}\end{equation}\end{Dn}
\begin{Dn}\label{Definition 5.5}Suppose that the wave operators $W_{\pm}(\mathcal{L}\mathcal{L}_{0})$
exist. They are complete if the relations
\begin{equation} Ran{W_{\pm}(\mathcal{L},\mathcal{L}_{0})}=\mathcal{H}_{ac}\label{5.6}\end{equation}
are valid.
\end{Dn}
Now we shall prove  the  central in this section  result:
\begin{Tm}\label{Theorem 5.6} Let the conditions of the  Theorem 5.2 be fulfilled. Then the corresponding wave operators $W_{\pm}(\mathcal{L},\mathcal{L}_{0})$ exist and are complete. \end{Tm}
We begin with the proof of the following lemma:
\begin{La}\label{Lemma 5.7} Let $g(r){\in}C^{\infty}_{0}$  then the equality
\begin{equation} e^{it\lambda}\int_{R^3}\phi^{*}_{n}(r,k,\lambda)g(r)dr=
\int_{R^3}\phi^{*}_{n}(r,k,\lambda)[e^{it\mathcal{L}}g(r)]dr,\quad
(1{\leq}n{\leq}4)\label{5.7}\end{equation}
is valid. Here $\lambda=-\sqrt{k^2+m^2}$, when $n=1,2,$ and $\lambda=\sqrt{k^2+m^2}$, when $n=3,4.$\end{La}
\emph{Proof.} We note that the vector function $\phi_{n}(r,k,\lambda)$ satisfies the equation \eqref{2.35}. Let us represent the left-hand  side of equality \eqref{5.7} as inner product in the Hilbert space $L^2(R^3)$:
\begin{equation}(e^{-it\lambda}\phi_{n}(r,k,\lambda),g(r))=
(e^{-it\mathcal{L}}\phi_{n}(r,k,\lambda),g(r)).\label{5.8}\end{equation}
The assertion of the lemma follows directly from \eqref{5.8}.\\
2. Let us write the partial case of relation \eqref{4.22}:
\begin{equation}\widetilde{g}_{0,n}(k,\lambda)=(2\pi)^{-3/2}
\int_{R^3}\phi_{0,n}^{*}(r,k,\lambda)g(r)dr,\quad (1{\leq}n{\leq}4)
\label{5.9}\end{equation}
where
\begin{equation}\phi_{0,n}(r,k,\lambda)=e^{ir{\cdot}k}\hat{g}_{n}(k),\quad
(1{\leq}n{\leq}4)\label{5.10}\end{equation}
and
$\lambda=-\sqrt{k^2+m^2}<-m,$ if  $n=1,2$, and  $\lambda=\sqrt{k^2+m^2}>m,$ if  $n=3,4$.
The vector functions $\phi_{0,n}(r,k,\lambda),\, (1{\leq}n{\leq}4)$ are solutions of the equation
\begin{equation}(\mathcal{L}_{0}-\lambda)\phi=0.\label{5.11}\end{equation}

\begin{La}\label{Lemma 5.8} Let the conditions of the  Theorem 5.2 be fulfilled.
If the $4{\times}1$ vector-function $g(r)$ belongs to
$L^{2}(R^3)$, then the  $4{\times}1$ vector-function $G(r)=W_{-}(\mathcal{L},\mathcal{L}_0)g(r)$ exists and \begin{equation} \hat{G}_{n}(k,\lambda)=\hat{g}_{0,n}(k,\lambda),\quad \lambda{\in}E_{+}.\label{5.12}\end{equation}
\end{La}
\emph{Proof.} Let us consider the expression $(f,W_{-}(\mathcal{L},\mathcal{L}_0)g)$.
We assume in addition that $g(r){\in}C^{\infty}_{0}$ and the support of $\hat{f}$ belongs to interval $(\alpha,\beta)$, where  $(\alpha,\beta){\in}E_{+}$. Using relation
\begin{equation}\frac{d}{dt}[e^{it\mathcal{L}}e^{-it\mathcal{L}_0}]=
ie^{it\mathcal{L}}Ve^{-it\mathcal{L}_0},\label{5.13}\end{equation}
we have
\begin{equation} (f,W_{-}(\mathcal{L},\mathcal{L}_0)g)-(f,g)=i\lim_{T{\to}-\infty}\int_{0}^{T}(f,e^{it\mathcal{L}}Ve^{-it\mathcal{L}_0}g)dt.
\label{5.14}\end{equation}
According to Abels limits we obtain  (see \cite{RS}, section 6, Lemma 5)
\begin{equation} (f,W_{-}(\mathcal{L},\mathcal{L}_0)g)-(f,g)=i\lim_{\delta{\to}+0}\int_{0}^{-\infty}e^{t\delta}(f,e^{it\mathcal{L}}Ve^{-it\mathcal{L}_0}g)dt.
\label{5.15}\end{equation}
We introduce
\begin{equation}F(r,t,\delta)=e^{t\delta}e^{it\mathcal{L}}Ve^{-it\mathcal{L}_0}g.
\label{5.16}\end{equation}
 Relations  \eqref{4.22} and \eqref{5.8} imply
\begin{equation}\hat{F}_{n}(k,\lambda,t,\delta)=(2\pi)^{-3/2}\int_{R^3}e^{t\delta}e^{it\lambda}\phi^{*}_{n}(r,k,\lambda)V(r)
[e^{-it\mathcal{L}_0}g(r)]dr, \label{5.17}\end{equation}
where $\lambda=-\sqrt{k^2+m^2}$, if $n=1,2$  and   $\lambda=\sqrt{k^2+m^2}$, if $ n=3,4.$
Now let us consider the integral (see \eqref{5.15}):
\begin{equation}J_{\delta}=\int_{0}^{-\infty}e^{t\delta}(f,e^{it\mathcal{L}}Ve^{-it\mathcal{L}_0}g)dt.
\label{5.18}\end{equation}
Taking into account  \eqref{4.50} and \eqref{5.17}  we have
\begin{equation}J_{\delta}=J_{\delta}(1)+J_{\delta}(2)+J_{\delta}(3)+J_{\delta}(4),\label{5.19}\end{equation}
where
\begin{equation}J_{\delta}(n)=\int_{0}^{-\infty}\int_{D}\hat{f}^{\star}_{n}(k)\hat{F}_{n}(k,\lambda,t,\delta)dkdt.
\label{5.20}\end{equation}
It follows from \eqref {3.3}, \eqref {3.5}, \eqref{5.17} and \eqref{5.20} that
\begin{equation}J_{\delta}(n)=-i\int_{D}\hat{f}^{\star}_{n}(k)(\phi_{n}(r,k,\lambda),V(r)
[B_{0}(\lambda-i\delta)g(r)])dk, \label{5.21}\end{equation}
 We receive from \eqref{5.21} and relation $B_{0}(\lambda-i\delta)=B_{0}^{*}(\lambda+i\delta)$ that
\begin{equation}lim J_{\delta}(n)=-i\int_{D}\hat{f}^{\star}_{n}(k)(B_{0}(\lambda)V(r)\phi_{n}(r,k,\lambda),
g(r))dk, \label{5.22}\end{equation}
where $\delta{\to}+0$.
Using \eqref{2.20} we rewrite the equality \eqref{5.22} in the form
\begin{equation}lim J_{\delta}(n)=-i\int_{D}\hat{f}^{\star}_{n}(k)([\phi_{0n}(r,k,\lambda)-\phi_{n}(r,k,\lambda)],
g(r)])dk, \label{5.23}\end{equation}
where $\delta{\to}+0$.
Hence we have \begin{equation}limJ_{\delta}(n)=-i\int_{D}\hat{f}^{\star}_{n}(k)[\hat{g}_{0n}(k)-\hat{g}_{n}(k)]
dk, \label{5.24}\end{equation}
where $\delta{\to}+0$. It follows from \eqref{5.15} and \eqref{5.24} that
\begin{equation}(f,W_{-}g)=\int_{D}[\hat{f}^{\star}_{3}(k)\hat{g}_{03}(k)+\hat{f}^{\star}_{4}(k)\hat{g}_{04}(k)]dk+\nonumber\end{equation}
\begin{equation}\int_{D}[\hat{f}^{\star}_{1}(k)\hat{g}_{01}(k)+\hat{f}^{\star}_{2}(k)\hat{g}_{02}(k)]dk
 \label{5.25}\end{equation}
The assertion of the
lemma follows directly from \eqref{5.25}.\\
\emph{Proof of the Theorem 5.6}\\ The wave operator $W_{-}(\mathcal{L},\mathcal{L}_{0})$
exists (see Theorem 5.2).
According to \eqref{4.49} and \eqref{5.12} we have $Ran[W_{-}(\mathcal{L},\mathcal{L}_{0})]=\mathcal{H}_{ac}$.
Hence, the  wave operator $W_{-}(\mathcal{L},\mathcal{L}_{0})$ is complete. In the same way we
can prove that the operator $W_{+}(\mathcal{L},\mathcal{L}_{0})$ is complete.
In this case, instead of operator $B_{+}(\lambda)$ (see \eqref{2.31}), we use the operator
\begin{equation}B_{-}(\lambda)f=\int_{R^{3}}V_{1}(r)B_{-}(r-s,\lambda)V_{1}(s)W_{1}(s)f(s)ds,
\label{5.26}\end{equation} where $B_{-}(r-s,\lambda)=B_{+}^{*}(s-r,\lambda)$. The theorem is proved.
\begin{Cy}\label{Corollary 5.9}Let the conditions of the Theorem 5.2 be fulfolled. Then the corresponding operator $S$ is unitary in the space $\mathcal{H}_{ac}$.\end{Cy}
\section{The scattering operator}

In this section we shall find the energetic representation $S(\lambda)$ of  the scattering  operator $S(\mathcal{L},\mathcal{L}_0).$ To do it we use  Definition 5.4
and  write
 \begin{equation} (f,(S-I)g)=((W_{+}-W_{-})f,W_{-}g).\label{6.1}\end{equation}
 We assume that conditions of the Theorem 5.2 are fulfilled. Hence, the wave operators $W_{\pm}$ and scattering operator $S$ exist. We assume also that $f(r)$ and $g(r)$ belong to $C_{0}^{\infty}$. Using  these assumptions we can change the order of integration in the next calculations.
It follows from  relation \eqref{6.1} and equality
\begin{equation}\frac{d}{dt}[e^{it\mathcal{L}}e^{-it\mathcal{L}_0}]=
ie^{it\mathcal{L}}Ve^{-it\mathcal{L}_0},\label{6.2}\end{equation}
that
\begin{equation}
(f,(S-I)g)=i\lim_{T{\to}\infty}\int_{-T}^{T}(e^{it\mathcal{L}}Ve^{-it\mathcal{L}_0}f,W_{-}g)dt.
\label{6.3}\end{equation}
According to Abels limits we obtain  (see \cite{RS}, section 6, Lemma 5)
\begin{equation} (f,(S-I)g)=i\lim_{\delta{\to}+0}\int_{-\infty}^{\infty}e^{-\delta|t|}
(e^{it\mathcal{L}}Ve^{-it\mathcal{L}_0}f,W_{-}g)dt.\label{6.4}\end{equation}
Let us introduce the vector functions
\begin{equation}F(r,t)=e^{it\mathcal{L}}Ve^{-it\mathcal{L}_0}f,\quad G(r)= W_{-}g.
\label{6.5}\end{equation}
Relations
\eqref{4.22} and  \eqref{5.25}  imply
\begin{equation}(F(r,t),G(r))= \int_{D}\sum_{n=1}^{4}[\hat{F}^{\star}_{n}(k,t)\hat{g}_{0n}(k)]dk.
 \label{6.6}\end{equation}
In view of \eqref{4.22}  and \eqref{6.5} we have
\begin{equation}\hat{F}_{n}(k,t)=(2\pi)^{-3/2}\int_{R^3}e^{it\lambda_{n}(k)}\phi^{*}_{n}(r,k,\lambda_{n}(k))V(r)
[e^{-it\mathcal{L}_0}f(r)]dr, \label{6.7}\end{equation}
where $\lambda_{n}(k)=-\sqrt{k^2+m^2},$ if $n=1,2$, and $\lambda_{n}(k)=\sqrt{k^2+m^2},$ if $n=3,4.$
It follows from \eqref{4.49} and Proposition 3.5 that
\begin{equation}f(r)=(2\pi)^{-3/2}\int_{D}e^{ir{\cdot}q}\sum_{s=1}^{4}[\hat{g}_{s}(q)\widetilde{f}_{s}(q)]dq.
\label{6.8}\end{equation}
Using \eqref{6.8} we have
\begin{equation}e^{-it\mathcal{L}_0}f(r)=(2\pi)^{-3/2}\int_{D}e^{ir{\cdot}q}\sum_{s=1}^{4}e^{-i\lambda_{s}(q)t}[\hat{g}_{s}(q)\widetilde{f}_{s}(q)]dq
 \label{6.9}\end{equation}
We need the following notations:
\begin{equation}T(q,k,\lambda)=(2\pi)^{-3}\int_{R^3}e^{-iqr}V(r)\Phi(r,k,\lambda)U_{0}^{*}(k)dr,\label{6.10}
\end{equation}
where
\begin{equation}U_{0}(k)=[\hat{g}_{1}(k),\hat{g}_{2}(k),\hat{g}_{3}(k),\hat{g}_{4}(k)],\label{6.11}\end{equation}
\begin{equation}\Phi(r,k,\lambda)=[\phi_{1}(r,k,\lambda),\phi_{2}(r,k,\lambda),\phi_{3}(r,k,\lambda),\phi_{4}(r,k,\lambda)].
\label{6.12}\end{equation}
Here
$\lambda=-\sqrt{k^2+m^2}$ if $n=1,2$ and $\lambda=\sqrt{k^2+m^2}$ if $n=3,4$. We note, that the vectors $\hat{g}_{s}(q)$ are orthogonal and $\|\hat{g}_{1}(q)\|=1.$
Hence the matrix  $U_{0}(q)$ is the unitary.
Integrating \eqref{6.4} over the variable t and using formulas  \eqref{6.5}-\eqref{6.12} we rewrite the
relation \eqref{6.4} in the form.
\begin{equation} (f,(S-I)g)=i\lim_{\delta{\to}+0}\sum_{1{\leq}s,n{\leq}4}J_{\delta}(s,n,\lambda).\label{6.13}
\end{equation}
Here
\begin{equation}J_{\delta}(s,n)=\nonumber\end{equation}
\begin{equation}\int_{D}[\int_{R^3}
\frac{2\delta}{\delta^2+[\lambda_{s}(q)-\lambda_{n}(k)]^2}
\overline{[F_{0}(q)]_{s}}T_{s,n}(q,k,\lambda)dq[G_{0}(k)]_{n}dk,\label{6.14}\end{equation}
where $T_{s,n}(q,k,\lambda)$  are elements of the matrix $T(q,k,\lambda),\,[F_{0}(q)]_{s}$  are elements
of the vector $U_{0}(q)\hat{f}(q),\,[G_{0}(k)]_{n}$  are elements of the vector $U_{0}(k)\hat{g}_{0}(k)$
and $\hat{g}_{0}(k)=col[\hat{g}_{01}(k),\hat{g}_{02}(k),\hat{g}_{03}(k)\hat{g}_{04}(k)].$
It is easy to see that
\begin{equation}\lim_{\delta{\to}+0}J_{\delta}(s,n,\lambda)=0, \label{6.15}\end{equation}
when $\lambda_{s}(k)=-\lambda_{n}(k).$ Let us consider the case when $\lambda_{s}(k)=\lambda_{n}(k).$
\begin{La}\label{Lemma 6.1} If $\lambda(|q|)=\sqrt{q^2+m^2}$, then
\begin{equation}Q=\lim_{\delta{\to}+0}\int_{0}^{\infty}
\frac{2\delta}{\delta^2+[\lambda(|q|)-\lambda(|k|)]^2}
d|q|=2{\pi}\sqrt{k^2+m^2}|k|^{-1}.\label{6.16}\end{equation}
\end{La}
\emph{Proof.}
We represent $\lambda(|q|)-\lambda(|k|)$ in the form
\begin{equation}\lambda(|q|)-\lambda(|k|)=\frac{q^2-k^2}{\lambda(|q|)+\lambda(|k|)}.\label{6.17}\end{equation}
The relation
\begin{equation}\lambda(|q|)-\lambda(|k|){\sim}\frac{|k|(|q|-|k|)}{\sqrt{k^2+m^2}}.\label{6.18}\end{equation}
is valid when $ |q|{\sim}|k|$ ,
Hence we receive the equality
\begin{equation}Q=(k^2+m^2)\lim_{\delta{\to}+0}\int_{0}^{\infty}
\frac{2\delta}{\delta^2(k^2+m^2)+k^2(|q|-|k|)^2}
d|q|.\label{6.19}\end{equation}
We use the formula
\begin{equation}\int\frac{dx}{ax^2+bx+c} =\frac{2}{\sqrt{4ac-b^2}}arctan\frac{2ax+b}{\sqrt{4ac-b^2}}.
\label{6.20}\end{equation}
Formula \eqref{6.20} is valid if $4ac-b^2>0$. In our case (see \eqref{6.19}) we have
$a=k^2,\,b=-2|k|^3,\,c=k^4+\delta^2(k^2+m^2),$ i.e.
\begin{equation} 4ac-b^2=4k^2\delta^2(k^2+m^2)>0.\label{6.21}\end{equation}
Relations  \eqref{6.20} and \eqref{6.21} imply
\begin{equation}Q=2{\pi}\sqrt{k^2+m^2}|k|^{-1}.\label{6.22}\end{equation}
The lemma is proved.\\
According to \eqref{6.10} the matrix function  $T_{s,n}(q,k)$ is bounded and continuous. Let us consider the case , when $\lambda_{s}(k)=\lambda_{n}(k).$ Then
using spherical system of coordinates in the space $q{\in}R^3$, Lemma 6.1 and relation \eqref{6.16} we obtain
\begin{equation}\lim_{\delta{\to}+0}J_{\delta}(s,n)=\nonumber\end{equation}
\begin{equation}\int_{D}a(|k|)
[\int_{S^2}\overline{[F_{0}(|k|\omega)]_{s}}T_{s,n}(|k|\omega,k,\lambda)d\Omega(\omega)[G_{0}(k)]_{n}dk,\label{6.23}\end{equation}
where $\omega=q/|q|$ and
\begin{equation}a(|k|)=2{\pi}\sqrt{k^2+m^2}|k|. \label{6.24}\end{equation}
Here by $S^2$ we denote the surface $|q|=1$ in the space $R^3$, $d\Omega$ is the  standard measure on the
surface $S^2$.
Let us introduce the $2{\times}2$ matrices
\begin{equation}T_{1}(q,k,\lambda)=\left(
                       \begin{array}{cc}
                         T_{1,1}(q,k,\lambda) & T_{1,2}(q,k,\lambda) \\
                         T_{2,1}(q,k,\lambda) & T_{2,2}(q,k,\lambda) \\
                       \end{array}
                     \right),\label{6.25}\end{equation}
where   $ \lambda=-\sqrt{k^2+m^2}.$
 \begin{equation}T_{2}(q,k,\lambda)=\left(
                       \begin{array}{cc}
                         T_{3,3}(q,k,\lambda) & T_{3,4}(q,k,\lambda) \\
                         T_{4,3}(q,k,\lambda) & T_{4,4}(q,k,\lambda) \\

                       \end{array}
                     \right),   \label{6.26}\end{equation}
where $ \lambda=\sqrt{k^2+m^2}.$ Using  \eqref{6.25} and \eqref{6.26} we construct the operators
\begin{equation}T_{1}(\lambda)h_{1}(\omega)=a(|k|)\int_{S^2}T_{1}(|k|\omega,|k|\omega^{\prime},\lambda)
{h}_{1}(\omega^{\prime})d\Omega(\omega^{\prime}),\label{6.27}\end{equation}
and
\begin{equation}T_{2}(\lambda)h_{2}(\omega)=a(|k|)\int_{S^2}T_{2}(|k|\omega,|k|\omega^{\prime},\lambda)
{h}_{2}(\omega^{\prime})d\Omega(\omega^{\prime}),\label{6.28}\end{equation}
\begin{equation}T_{s,n}(\lambda)h(\omega)=a(|k|)\int_{S^2}T_{s,n}(|k|\omega,|k|\omega^{\prime},\lambda)
{h}(\omega^{\prime})d\Omega(\omega^{\prime}),\label{6.29}\end{equation}
where    $h_{1}(\omega)$ and   $h_{2}(\omega)$    are   $2{\times}1$ vector functions, $h(\omega)$ is a function, $\omega=q/|q|,\, \omega^{\prime}=k/|k|.$
Let us  introduce the Hilbert space $\mathcal{M}$ of the vector functions $h(\omega)=col[h_{1}(\omega),h_{2}(\omega)]$. The norm in the space $\mathcal{M}$ is defined by the relation
\begin{equation}\|h\|^{2}=\int_{S^2}(|h_{1}(\omega)|^{2}+|h_{2}(\omega)|^{2})d\Omega(\omega).\label{6.30}\end{equation}
The subset of $\mathcal{M}$, which corresponds to absolutely continuous spectrum,  we denote by $\mathcal{M}_{ac}.$ The operators $T_{1}(\lambda)$ and $T_{2}(\lambda)$  act in the space $\mathcal{M}_{ac}.$
 Relations \eqref{6.13}, \eqref{6.15} and \eqref{6.23}-\eqref{6.29} imply the assertion
\begin{Pn}\label{Proposition 6.2}Let conditions of the Theorem 5.2 be fulfilled. The energetic representation $S(\lambda)$ of  the scattering  operator $S(\mathcal{L},\mathcal{L}_0)$ has the form
\begin{equation}S(\lambda)= S_{1}(\lambda),\, \lambda<-m; \quad S(\lambda)= S_{2}(\lambda),\, \lambda>m,
\label{6.31}\end{equation}
where
\begin{equation}S_{1}(\lambda)=I+iT_{1}(\lambda),\quad S_{2}(\lambda)=I+iT_{2}(\lambda).\label{6.32}\end{equation}
\end{Pn}
Now let us investigate   the structure of the operators
 $T_{s,n}(\lambda)$. To do it we introduce the operator
\begin{equation}\mathcal{F}_{s}(\lambda)f(r)=\sqrt{a(|k|)}\int_{R^3}
e^{-iq{\cdot}r}\hat{g}_{s}^{*}(q)V_1(r)f(r)dr.
\label{6.33}\end{equation}
The adjoint to $\mathcal{F}_{s}(\lambda)$ operator has the form
\begin{equation}\mathcal{F}^{\star}_{s}(\lambda)h(\omega)=\sqrt{a(|k|)}V_{1}(r)\hat{g}_{s}(r)\int_{S^2}e^{i|k|\omega{\cdot}r}
h(\omega)d\Omega(\omega)
\label{6.34}\end{equation}
The following assertion is valid.
\begin{Pn}\label{Proposition 6.3}Let conditions of Theorem 5.2 be fulfilled. If $\lambda$  belong to
$E_{+}$ then  the operator $T_{s,n}(\lambda)$ can be represented in the form
\begin{equation}T_{s,n}(\lambda)=\mathcal{F}_{s}(\lambda)W_{1}(r)[I+(2\pi)^{3/2}B_{+}(\lambda)]^{-1}
\mathcal{F}^{\star}_{n}(\lambda)\label{6.35}\end{equation}\end{Pn}
\emph{Proof.}
 According to \eqref{2.30} we have
\begin{equation}\psi_{n}(r,k)=[I+(2\pi)^{3/2}B_{+}(\lambda)]^{-1}e^{ikr}V_{1}(r)\hat{g}_{n}(k)
\label{6.36}
\end{equation}Relation \eqref{6.35} follows directly from
\eqref{6.32}-\eqref{6.34} and \eqref{6.36}.

\begin{Tm}\label{Theorem 6.4}Let conditions of Theorem 5.2 be fulfilled. Then for each
$\lambda{\in}E_{+}$  the operators $T_{1}(\lambda)$ and $T_{2}(\lambda)$ belong to the Hilbetrt- Schmidt class.\end{Tm}
\emph{Proof.} We use the relation
\begin{equation}[\mathcal{F}^{\star}_{n}\mathcal{F}_{n}f]=\nonumber\end{equation}
\begin{equation}{4\pi} a(|k|)\int_{R^3}V_{1}(r)\hat{g}_{n}(|k|)
\frac{sin(|k||r-s|)}{|k||r-s|}\hat{g}_{n}^{\star}(|k|)V_{1}(s)
f(s)ds.\label{6.37}\end{equation}
The  function $\|V(r)\|$ belongs to Rollnik class (see \cite{Sim},Ch.1). Hence  it follows from \eqref {6.36} that the operators
$\mathcal{F}^{\star}_{n}\mathcal{F}_{n}$ belong to the Hilbert-Schmidt class.
\begin{Cy}\label{Corollary 6.5}Let conditions of the Theorem 5.2 be fulfilled. Then for each $\lambda{\in}E_{+}$ the eigenvalues
$\mu_{j}(\lambda)$ of $S(\lambda)$  are such that
\begin{equation}\sum_{j=1}^{\infty}[\mu_{j}(\lambda)-1]^2<\infty .
\label{6.38}\end{equation}\end{Cy}
\emph{Proof.} The operators
$T_{1}(\lambda)$  and $T_{2}(\lambda)$ belong  to the Hilbert-Schmidt class. Hence, the relation \eqref{6.38} is valid.\\

\section{Scattering amplitude}
In this section we shall investigate  the connection between the stationary scattering problem
(scattering amplitude)   and dynamical scattering problem (scattering operator).\\
\emph{We suppose that the conditions of Theorem 2.2 and condition \eqref{2.32} are fulfilled.}\\
Let us  introduce the Hilbert space $\mathcal{H}$ of the vector functions $h(\omega)=col[h_{1}(\omega),h_{2}(\omega)]$. The norm in the space $\mathcal{H}$ is defined by the relation
\begin{equation}\|h\|^{2}=\int_{S^2}(|h_{1}(\omega)|^{2}+|h_{2}(\omega)|^{2})d\Omega(\omega).\label{7.1}\end{equation}
We consider the functions
\begin{equation}f_{s,n}(q,k)=\hat{g}_{s}^{*}(q)f(\omega,k,n),\label{7.2}\end{equation}
where the vector function $f(\omega,k,n)$ is defined by the relation \eqref{2.34}.
Now we construct the following $2{\times}2$ matrices
\begin{equation}F_{1}(q,k)=\left(
                       \begin{array}{cc}
                         f_{1,1}(q,k) & f_{1,2}(q,k) \\
                         f_{2,1}(q,k) & f_{2,2}(q,k) \\
                       \end{array}
                     \right),\label{7.3}\end{equation}
 \begin{equation}F_{2}(q,k)=\left(
                       \begin{array}{cc}
                         f_{3,3}(q,k) & f_{3,4}(q,k) \\
                         f_{4,3}(q,k) & f_{4,4}(q,k) \\

                       \end{array}
                     \right)   \label{7.4}\end{equation}
Let us compare formulas \eqref{2.34}, \eqref{6.10}  and \eqref{7.2}. It is easy to see that
\begin{equation}
F_{p}(q,k)=-\lambda{\gamma}T_{p}(q,k),
\quad \gamma=2^{4}\pi^{5},\quad p=1,2.\label{7.5}\end{equation}
According to Corollary 5.9 the operators $S_{1}(\lambda)$ and $S_{2}(\lambda)$ are unitary in the space $\mathcal{H}_{ac}$. Hence there exists a complete orthonormal systems of eigenvectors $G_{j,1}(\omega,\lambda)$ and $G_{j,2}(\omega,\lambda)$  of the operators $S_{1}(\lambda)$ and $S_{2}(\lambda)$ respectively. We denote the corresponding eigenvalues by $\mu_{j,1}(\lambda)$ and  $\mu_{j,2}(\lambda)$. We note that  $|\mu_{j,1}(\lambda)|=|\mu_{j,2}(\lambda)|=1.$
The vector functions
\begin{equation}F_{1,1}=[f_{1,1},f_{1,2}],\quad F_{2,1}=[f_{2,1},f_{2,2}],
 \label{7.6}\end{equation}
and the vector functions
\begin{equation}F_{1,2}=[f_{3.3},f_{3.4}],\quad F_{2,2}=[f_{4.3},f_{4.4}]
 \label{7.7}\end{equation}
 belong to the space  $\mathcal{H}$. Further we consider the case when
\begin{equation}\mathcal{H}_{ac}=\mathcal{H}.\label{7.8}\end{equation}  Hence  the vector functions $F_{m,p}$ can be represented in the form of series:
 \begin{equation}F_{m,p}(\omega,\omega^{\prime},\lambda)=\sum_{j}a_{j}(m,p,\omega,\lambda)G_{j,p}^{*}(\omega^{\prime},\lambda),
 \label{7.9}\end{equation}
where $q=|k|\omega$, $k=|k|\omega^{\prime}$ and
\begin{equation}a_{j}(m,p,\omega,\lambda)=\int_{S^2}F_{m,p}(\omega,\omega^{\prime},\lambda)G_{j,p}(\omega^{\prime},\lambda)
d\Omega(\omega^{\prime}) \label{7.10}\end{equation}

It follows from \eqref{6.28} and \eqref{7.5}-\eqref{7.10} the assertion.
\begin{Tm}\label{Theorem 7.1}Let the conditions of Theorem 2.2 and conditions \eqref{2.32}, \eqref{7.8} be fulfilled.
Then we have
\begin{equation}
F_{p}(\omega,\omega^{\prime},\lambda)=i\lambda{\gamma}\sum_{j}(\mu_{j,p}(\lambda)-1)G_{j.p}(\omega,\lambda)G_{j,p}^{*}(\omega^{\prime},\lambda),
 \quad p=1,2.\label{7.11}\end{equation}\end{Tm}
Let us define the total cross sections
\begin{equation}\sigma_{p}(\lambda)=\int_{S^2}\int_{S^2}F_{p}(\omega,\omega^{\prime},\lambda)
F_{p}^{*}(\omega,\omega^{\prime},\lambda)d\Omega(\omega)d\Omega(\omega^{\prime}),\quad p=1,2.
\label{7.12}
\end{equation}
\begin{Cy}\label{Corollary 7.2}Let conditions of Theorem 2.2 and conditions \eqref{2.32}, \eqref{7.8} be fulfilled. Then the total cross sections have  the forms
\begin{equation}\sigma_{p}(\lambda)=(\lambda\gamma)^2\sum_{j}|\mu_{j}(\lambda)-1|^{2}\int_{S^2}
G_{j.p}(\omega,\lambda)G_{j.p}^{*}(\omega,\lambda)d\Omega(\omega).\label{7.13}\end{equation}\end{Cy}
\begin{Rk}\label{Remark 7.3} The operators $T_{p}(\lambda)$  belong to the Hilbert-Schmidt class and
\begin{equation}
Tr[\int_{S^2}G_{j.p}(\omega,\lambda)G_{j.p}^{*}(\omega,\lambda)d\Omega(\omega)]=1.
\label{7.14}\end{equation}
 Hence  the series \eqref{7.11} and \eqref{7.13} converge.\end{Rk}
 It follows from \eqref{7.13} and  \eqref{7.14} that
\begin{equation}Tr[\sigma_{p}(\lambda)]=(\lambda\gamma)^2\sum_{j}|\mu_{j,p}(\lambda)-1|^{2}<\infty
\label{7.15}\end{equation}
\begin{Rk}\label{Remark 7.4}Formulas \eqref{7.11}, \eqref{7.13} and \eqref{7.15} give the connections between the stationary scattering results  $(F_{p}(\omega,\omega^{\prime},\lambda))$ and the dynamical scattering results  $(G_{j,p}(\omega,\lambda),\,\mu_{j,p}(\lambda))$. So, formulas \eqref{7.11} and \eqref{7.13}are quantum mechanical analogues of the ergodic formulas in classical mechanics.Vor radial case , when $V(r)=V(|r|)$, the
quantum mechanical analogues of the ergodic formulas  were obtained in our papers  \cite{Sakh7},\cite{Sakh8}.\end{Rk}

\section{Schr\"odinger operators}
1. Let us consider the Schr\"odinger operators
\begin{equation}\mathcal{L}u(r)=-\Delta{u}+V(r)u,\quad \mathcal{L}_{0}u(r)=-\Delta{u},\label{8.1}\end{equation}
where $r=(r_1,r_2,r_3){\in}R^3$. In the present  section we shall
prove that  formulas of the type
\eqref{7.11} and \eqref{7.15} are valid for the Schr\"odinger operators too.
Let us write the Lippmann- Schwinger equation:
\begin{equation}\phi(r,k)=e^{ik{\cdot}r}-\frac{1}{4\pi}\int_{R^{3}}\frac{e^{i|k||r-s|}}{|r-s|}V(s)\phi(s,k)ds,
\label{8.2}\end{equation}
The modified Lippmann-Schwinger equation has the form
\begin{equation}(I+K(\lambda))\psi(r,k)=e^{ik{\cdot}r}|V(r)|^{1/2},\label{8.3}\end{equation}
where the operator $K(\lambda)$ is defined by the relation
\begin{equation}K(\lambda)f(r)=\frac{1}{4\pi}\int_{R^{3}}|V(r)|^{1/2}\frac{e^{i|k||r-s|}}{|r-s|}W(s)|V(s)|^{1/2}f(s)ds.
\label{8.4}\end{equation}Here $W(s)=sgn{V(s)}.$
\begin{Rk} \label{Remark 8.1} Further we assume that the potential $V(r)$ belongs to the Rollnik class,
$V(r)=\overline{V(r)}$ and $V(r){\in}L(R^3)$.\end{Rk}
\begin{Dn}\label{Definition 8.2}The point $\lambda>0$ is an exceptional value if the equation $[I+K(\lambda)]\psi=0$ has nontrivial solution in the space $L^2(R^3)$.\end{Dn}
 We denote by $\mathcal{E}_{+}$ the set of exceptional points and we denote by $E_{+}$ the set
 of such points $\lambda>0$ that  $\lambda{\notin}\mathcal{E}_{+}.$\\
 \begin{La}\label{Lemma 8.3} If $\lambda{\in}E_{+}$, then equation \eqref{8.3} has one and only one
solution $\psi(r,k)$ in $L^2(R^3)$. \end{La}
\begin{Cy}\label{Corollary 8.4} If $\lambda{\in}E_{+}$
, then equation \eqref{8.2} has one and only one
solution $\phi(r,k)$ which satisfies the condition $|V(r)|^{1/2}\phi(r,k){\in}L^2(R^3)$. \end{Cy}
We formulate the following result (see \cite{RS}, p.115).
 \begin{Cy}\label{Corollary 8.5}
1) The operator $\mathcal{L}$ has only discrete spectrum in the domain $\lambda<0$.\\
2) The discrete spectrum  of the  operator $\mathcal{L}$ has no limit points in the  $\lambda<0$.\\
3)The set $\mathcal{E}_{+}$ is  bounded,  closed and has  Lebesgue  measure equal to zero. \end{Cy}
Let us consider the scattering operator $S(\lambda)$ in the energetic representation. It is
proved (see \cite {RS}, p.110) that $S(\lambda)$ is unitary. Let us introduce the operator function
\begin{equation}T(\lambda)=(2i\pi)^{-1}[I-S(\lambda)]. \label{8.5}\end{equation}
We need the following assertion (\cite{Kato})
\begin{Pn}\label{Proposition 8.6} For each
$\lambda{\in}E_{+}$  the operator $T(\lambda)$  belongs to the Hilbert-Schmidt class and  can be represented in the form
\begin{equation}T(\lambda)=\mathcal{F}(\lambda)W(r)[I+K(\lambda)]^{-1}
\mathcal{F}^{\star}(\lambda)\label{8.6}\end{equation}\end{Pn}
 The operator $\mathcal{F}(\lambda)$ is defined by the formula
\begin{equation}\mathcal{F}(\lambda)f(r)=\mu\int_{R^3}
e^{-i|k|\omega{\cdot}r}|V(r)|^{1/2}f(r)dr,\quad |k|^2=\lambda,
\label{8.7}\end{equation}
where $\mu=(1/4)(\lambda)^{1/4}\pi^{-3/2}$.
The adjoint to $\mathcal{F}(\lambda)$ operator has the form
\begin{equation}\mathcal{F}^{\star}(\lambda)h(\omega)=\mu\int_{S^2}e^{i|k|\omega{\cdot}r}|V(r)|^{1/2}h(\omega)
d\Omega(\omega)
\label{8.8}\end{equation}Here by $S^2$ we denote the surface $|\omega^{\prime}|=1$ in the space $R^3$, $d\Omega$ is the  standard measure on the
surface $S^2$.
It follows from \eqref{8.6} that the operator $T(\lambda)$ can be represented in the form
\begin{equation}T(\lambda)h(\omega^{\prime})={\mu^2}\int_{S^2}T(\omega,\omega^{\prime},\lambda)h(\omega^{\prime})
d\Omega(\omega^{\prime}),
  \label{8.9}\end{equation}
where $ \omega=k/|k|,\, \omega^{\prime}=k^{\prime}/|k^{\prime}|,\,|k|=|k^{\prime}|=\sqrt{\lambda}.$
Changing the order of integrals in \eqref{8.6} and using \eqref{8.3} we obtain
\begin{equation} T(\omega,\omega^{\prime},\lambda)=\int_{R^3}e^{-i\sqrt{\lambda}(s{\cdot}\omega)}
V(s)\phi(s,k^{\prime})ds,\label{8.10}\end{equation}
Further we need the following result (see \cite{IK}, p.11).
\begin{Tm}\label{Theorem 8.7}Let the function $V(r)$ belong to space $L^{2}(R^3)$.
 If
\begin{equation}|V(r)|=O(|r|^{-3-\delta}),\quad \delta>0,\quad |r|{\to}\infty,\label{8.11}\end{equation}
then the solution $\phi(r,k)$ of LS equation \eqref{8.2} has the form
\begin{equation}\phi(r,k)=e^{ik{\cdot}r}+
\frac{e^{i\sqrt{\lambda}|r|}}{|r|}f(\omega,\omega^{\prime},\lambda)
+o(1/|r|),\quad |r|{\to}\infty , \label{8.12}\end{equation}
where $\omega=r/|r|,\,\omega^{\prime}=k/|k|$ and
\begin{equation}f(\omega,\omega^{\prime},\lambda)=-\frac{1}{4\pi}\int_{R^3}e^{-i\sqrt{\lambda}(s{\cdot}\omega)}
V(s)\phi(s,k)ds.  \label{8.13}\end{equation}
\end{Tm}
The function $f(\omega,\omega^{\prime},\lambda)$ is  the scattering amplitude.\\
We note that $\omega$ and $\omega^{\prime}$ are defined in the cases \eqref{8.9} and \eqref{8.12}
differently.
It is easy to see that
\begin{equation}
f(\omega,\omega^{\prime},\lambda)=-\frac{1}{4\pi}T(\omega,\omega^{\prime},\lambda),
\label{8.14}\end{equation}
2. We introduce the Hilbert space $\mathcal{H}$ of the functions $h(\omega)$. The norm in the space $\mathcal{H}$ is defined by the relation
\begin{equation}\|h\|^{2}=\int_{S^2}|h(\omega)|^{2}d\Omega(\omega).\label{8.15}\end{equation}
Let us compare formulas \eqref{8.6} and \eqref{8.3}, \eqref{8.11}.
We note that the operator $S(\lambda)$ is unitary in the space $\mathcal{H}$. Hence there exists a complete orthonormal system of eigenfunctions $G_{j}(\omega,\lambda)$ of the operator $S(\lambda)$. We denote the corresponding eigenvalues by $\mu_{j}(\lambda)$, where $|\mu_{j}(\lambda)|=1.$
The  function $f(\omega,\omega^{\prime},\lambda)$, where $\omega$ and $\lambda$ are fixed,
 belongs to the space  $\mathcal{H}$. Hence  the  function $f(\omega,\omega^{\prime},\lambda)$ can be represented in the form of series:
 \begin{equation}f(\omega,\omega^{\prime},\lambda)=\sum_{j}a_{j}(\omega)\overline{G_{j}(\omega^{\prime},\lambda)},
 \label{8.16}\end{equation}
where
\begin{equation}a_{j}(\omega)=\int_{S^2}f(\omega,\omega^{\prime},\lambda)G_{j}(\omega^{\prime},\lambda)
d\Omega(\omega^{\prime}) \label{8.17}\end{equation}
It follows from \eqref{8.5},\eqref{8.9} and \eqref{8.15}-\eqref{8.17}  the assertion.
\begin{Tm}\label{Theorem 8.8}Let conditions of Theorem 8.7 be fulfilled. Then we have
\begin{equation}f(\omega,\omega^{\prime},\lambda)=\frac{2\pi}{i\sqrt{\lambda}}\sum_{j}(\mu_{j}(\lambda)-1)G_{j}(\omega,\lambda)
\overline{G_{j}(\omega^{\prime},\lambda)},
\label{8.18}\end{equation}\end{Tm}
Let us define the total cross section
\begin{equation}\sigma(\lambda)=\int_{S^2}\int_{S^2}|f(\omega,\omega^{\prime},\lambda)|^{2}d\Omega(\omega)d\Omega(\omega^{\prime}).
\label{8.19}\end{equation}
\begin{Cy}\label{Corollary 8.9}Let conditions of Theorem 8.7 be fulfilled.Then the total cross section has the form
\begin{equation}\sigma(\lambda)=\frac{4\pi^{2}}{\lambda}\sum_{j}|\mu_{j}(\lambda)-1|^{2}.\label{8.20}\end{equation}\end{Cy}
\begin{Rk}\label{Remark 8.10} The operator $T(\lambda)$  belongs to the Hilbert-Schmidt class. Hence  the series \eqref{8.18} and \eqref{8.20} converge.\end{Rk}
\begin{Rk}\label{Remark 8.11}The results of type \eqref{8.18} and \eqref{8.20} are well-known for the radial case , when $V(r)=V(|r|)$.\end{Rk}
\begin{Rk}\label{Remark 8.12}Formulas \eqref{8.18} and \eqref{8.20} give the connections between the stationary scattering results  $(f(\omega,\omega^{\prime},\lambda))$ and the dynamical scattering results  $(G_{j}(\omega,\lambda),\,\mu_{j}(\lambda))$. So, formulas \eqref{8.18} and \eqref{8.20}are quantum mechanical analogues of the ergodic formulas in classical mechanics.Vor radial case , when $V(r)=V(|r|)$, the
quantum mechanical analogues of the ergodic formulas  were obtained in our papers  \cite{Sakh7},\cite{Sakh8}.\end{Rk}

\end{document}